\documentclass[11pt]{article}
\usepackage{array, amsmath, amssymb, latexsym}

\begin{document}
\begin{center}
{\Large\bf  Combinatorial study on the group of parity alternating \\ \vspace{0.5cm}
permutations }  \\
\vspace{0.9cm}
{\bf Shinji Tanimoto}\\ 
\vspace{0.5cm}
Department of Mathematics,
University of Kochi,\\
Kochi 780-8515, Japan\footnote{Former affiliation}.  \\
\end{center}
\begin{abstract}
The combinatorial theory for the set of parity alternating permutations is expounded. In view of the numbers of
ascents and inversions, several enumerative aspects of the set are investigated. In particular, it is shown
that signed Eulerian numbers have intimate relationships to the set.
\end{abstract}
\vspace{0.7cm}
{\large \bf 1. Introduction} \\  
\\
\indent
Let $n$ be a positive integer. A permutation of the set $[n]= \{1, 2, \ldots, n\}$ will be called a {\it parity alternating permutation} ({\it PAP}),
if its entries assume even and odd integers alternately, such as $72581634$ or $1452763$. 
It is readily checked that such permutations form a subgroup of the symmetric group, ${\mathcal S}_n$, of degree $n$ 
and it will be denoted by ${\mathcal P}_n$. The objective of this paper is to study combinatorial
properties of permutations in ${\mathcal P}_n$ with respect to permutation statistics. \\
\indent
The extreme ends $a_1$ and $a_n$ in a PAP $\xi = a_{1}a_{2}\cdots a_{n}$ are both odd, when $n$ is odd. 
Hence the cardinality $|{\mathcal P}_n|$ is equal to
\[
            \left(\frac{n+1}{2}\right)!\left(\frac{n-1}{2}\right)! =
       \frac{n+1}{2} \left(\left(\frac{n-1}{2}\right)!\right)^2, 
\] 
when $n$ is odd, while it is equal to 
\[
                    2\left(\left(\frac{n}{2}\right)!\right)^2,
\]
when $n$ is even, because there are two possibilities; either even $a_1$ or odd $a_1$. \\
\indent
Let $\xi = a_1a_2\cdots a_n$ be a permutation of $[n]$.
An {\it ascent} of $\xi$ is an adjacent pair such that $a_i < a_{i+1}$ 
for some $i$ ($1\le i \le n-1$).  
Let ${\mathcal S}_{n,k}$ be the set of all permutations of $[n]$ with exactly $k$ ascents, where $0 \le k \le n-1$. 
Its cardinality, $|{\mathcal S}_{n,k}|$, is the classical {\it Eulerian number} (see [1] or [3]).
In relation to it let us define the set of PAPs with $k$ ascents by
\[
{\mathcal P}_{n, k} = {\mathcal P}_n \cap {\mathcal S}_{n,k}.
\]
\indent
For our study of ${\mathcal P}_{n, k}$ we will make use of an operator introduced in [4]. 
It is useful for a study of permutations in view of the ascent number (see [4--7]). 
The operator $\sigma$ is defined by 
adding one to all entries of a permutation $a_{1}a_{2}\cdots a_{n}$ of $[n]$, but by changing $n+1$ into one. 
However, when $n$ appears at either end of a permutation, it is removed
and one is put at the other end. To be more precise, we define as follows.
{\begin{itemize}
\item[\rm(a)] 
$\sigma (a_1a_{2}\cdots a_{n}) =b_1b_2 \cdots b_n$, 
where $b_{i}= a_{i}+1$ for $1\le i \le n$ and $n+1$ is replaced by one at the position, {\it i.e.}, $b_i = 1$, when 
$a_{i} = n$ for $2\le i \le n-1$. 
\item[\rm (b)] 
$\sigma (a_{1}\cdots a_{n-1}n) = 1 b_1b_2\cdots b_{n-1}$, 
where $b_{i}=a_{i}+1$ for $1\le i \le n-1$. 
\item[\rm (c)]
$\sigma (na_1\cdots a_{n-1}) = b_1b_2\cdots b_{n-1}1$, where $b_{i}=a_{i}+1$
for $1\le i \le n-1$. 
\end{itemize}}
\indent
This operator is a bijection on ${\mathcal S}_n$ and preserves the number of ascents of permutations, as is easily observed. 
Hence it becomes a bijection on ${\mathcal S}_{n, k}$ onto itself for each ascent number $k$. For a permutation 
$\xi$ and a positive integer $\ell$, we denote successive applications of $\sigma$ by 
$\sigma^\ell \xi=\sigma(\sigma ^{\ell -1}\xi)$, $\sigma^0$ being the identity operator on permutations. \\
\indent
It is easy to see that in the case of even $n$ the operator $\sigma$ is also a bijection on ${\mathcal P}_{n,k}$ onto itself.
When $n$ is odd, however, this is not the case. Taking a PAP $\xi = 1476523$, for example, 
$\sigma \xi = 2517634$ is not a PAP. 
The following table shows the cardinalities $|{\mathcal P}_{n, k}|$ for small $n$, the top row being the ascent number $k$. 
\\ 
\begin{center}
\begin{tabular}{l|rrrrrrrrrrr}
         \noalign{\hrule height 0.8pt}
        $|{\mathcal P}_{n, k}|$ & 0 & 1 & 2 & 3 & 4 & 5 & 6 & 7 & 8 & ~9 \\
       \hline
         $n = 2$ & 1 & 1 \\
         $n = 3$ & 1 & 0 & 1  \\
         $n = 4$ & 1 & 3 & 3 & 1 \\
         $n = 5$ & 1 & 2 & 6 & 2 & 1 \\
         $n = 6$ & 1 & 9 & 26 & 26 & 9 & 1 \\
         $n = 7$ & 1 & 8 & 39 & 48 & 39 & 8 & 1 \\ 
         $n = 8$ & 1 & 23 & 165 & 387 & 387 & 165 & 23 & 1 \\ 
         $n= 9$ & 1 & 22 & 228 & 674 & 1030 & 674 & 228 & 22 &~1 \\
         $n= 10$ & 1 & 53 & 860 & 4292 & 9194 & 9194 & 4292 & 860 & 53 &~1 \\
	 \noalign{\hrule height 0.8pt}
\end{tabular}
\end{center}
\vspace{0.6cm}
\indent
\indent
For permutations $\xi = a_1a_2\cdots a_n$ we can define an involution operator by $\xi^{\ast} = a_n \cdots a_2 a_1$,
{\it i.e.}, ${\xi^{\ast}}^{\ast} = \xi$. If $\xi$ has $k$ ascents, then $\xi^{\ast}$ has $n - 1 -k$ ascents, for 
$\xi$ has $n-1$ adjacent pairs.  \\
\indent
An {\it inversion} of a permutation $\xi = a_1a_2\cdots a_n$ is a 
pair $(i, j)$ such that $1 \le i < j \le n$ and $a_i > a_j$. 
Let us denote by ${\rm inv}(\xi)$ the number of inversions in $\xi$. 
A permutation is called even or odd if it has an even or odd number of inversions, respectively. 
We observe the number of inversions of a permutation when $\sigma$ is applied. 
When $n$ appears at either end of a permutation $\xi$ as in (b) or (c), it is evident that 
\[
{\rm inv}(\sigma \xi) = {\rm inv}(\xi).
\]
As for the case (a), i.e., $a_{i}=n$ for some $i$ ($2\le i \le n-1$), we get
$\sigma (a_1a_{2}\cdots a_{n}) =b_1b_2 \cdots b_n$ and $b_{i}= 1$.  In this case, $n-i$ inversions 
$(i, i+1), \ldots, (i, n)$ of $\xi$ vanish and, in turn, 
$i-1$ inversions $(1, i), \ldots, (i-1, i)$ of $\sigma \xi$ occur. Hence
the difference between the numbers of inversions is given by
\begin{eqnarray}
 {\rm inv}(\sigma \xi) - {\rm inv}(\xi) = (i-1) - (n-i) = 2i - (n + 1).
\end{eqnarray}
This means that, when $n$ is even, each application of the operator $\sigma$ changes the parity 
of permutations as long as $n$ remains in the interior of permutations. 
If $n$ is odd, however, the operator also preserves the parity of all permutations of $[n]$.\\
\indent
A {\it signed Eulerian number} is the difference of the numbers of even permutations and 
odd ones in ${\mathcal S}_{n,k}$, which is denoted by $d_{n,k}$. Namely, letting
${\mathcal S}^{\rm e}_{n,k}$ and ${\mathcal S}^{\rm o}_{n,k}$ be the sets of even permutations 
and odd ones in ${\mathcal S}_{n,k}$, respectively, it is defined by 
\[  d_{n,k} = |{\mathcal S}^{\rm e}_{n,k}| - |{\mathcal S}^{\rm o}_{n,k}|.
\]
As was shown in [2], the recurrence relation for signed Eulerian numbers is given by
\begin{eqnarray}
     d_{n,k} = 
        \left\{\begin{array}{ll}
                    (n-k)d_{n-1,k-1}+ (k+1)d_{n-1,k}, & \mbox{if  $n$ is odd}, \\
                    d_{n-1,k-1} - d_{n-1,k}, &  \mbox{if  $n$ is even},\\
                    \end{array}  \right.
\end{eqnarray}
with initial condition $d_{1,0}=1$. \\
\\
\\
\\
{\large \bf 2.  Enumerative properties of ${\mathcal P}_n$ for even $n$}  \\
\\
\indent
In this section we assume that $n$ is an even positive integer and
we derive several enumerative properties of ${\mathcal P}_{n, k}$ of parity alternating permutations with $k$ ascents. 
We subdivide the set ${\mathcal P}_{n, k}$ into the following four classes. 
We denote by ${\mathcal E}^{\uparrow}_{n, k}$ (or ${\mathcal E}^{\downarrow}_{n, k}$)
the set of permutations $a_{1}a_{2}\cdots a_{n} \in {\mathcal P}_{n, k}$ such that
$a_1$ is even and $a_n > a_1$ (or $a_n < a_1$). 
Similarly we denote by ${\mathcal O}^{\uparrow}_{n, k}$ (or ${\mathcal O}^{\downarrow}_{n, k}$)
the set of permutations $a_{1}a_{2}\cdots a_{n} \in {\mathcal P}_{n, k}$ such that
$a_1$ is odd and $a_n > a_1$ (or $a_n < a_1$).  \\
\indent
According to (a)--(c) it is easily shown that $\sigma$ maps ${\mathcal P}^{\uparrow}_{n, k} = {\mathcal O}^{\uparrow}_{n, k} \cup 
{\mathcal E}^{\uparrow}_{n, k}$ 
onto ${\mathcal P}^{\uparrow}_{n, k}$. It also maps ${\mathcal P}^{\downarrow}_{n, k} = {\mathcal O}^{\downarrow}_{n, k} 
\cup {\mathcal E}^{\downarrow}_{n, k}$ 
onto ${\mathcal P}^{\downarrow}_{n, k}$.  Note that ${\mathcal P}_{n, k} = {\mathcal P}^{\uparrow}_{n, k} 
\cup {\mathcal P}^{\downarrow}_{n, k}$, whose cardinality is denoted by $|{\mathcal P}_{n, k}|$. \\
\indent
As was shown in [4], to each permutation $\xi$ of $[n]$ there corresponds a smallest positive 
integer $\ell$ such that $\sigma^{\ell}\xi = \xi$, which is the {\it period} of $\xi$ and denoted by $\pi(\xi)$. Its trace 
\[
           \{\sigma \xi, \sigma^2 \xi, \ldots, \sigma^{\pi(\xi)}\xi=\xi \} 
\]
is called the {\it orbit} of $\xi$. \\
\indent
For an even $n$, $\xi \in {\mathcal P}^{\uparrow}_{n, k}$ holds if and only if $\sigma \xi \in {\mathcal P}^{\uparrow}_{n, k}$, and 
$\xi \in {\mathcal P}^{\downarrow}_{n, k}$ holds if and only if $\sigma \xi \in {\mathcal P}^{\downarrow}_{n, k}$. Therefore, we can
consider orbits both in ${\mathcal P}^{\uparrow}_{n, k}$ and in ${\mathcal P}^{\downarrow}_{n, k}$ under $\sigma$.
Permutations of the form $a_1a_2 \cdots a_{n-1}n \in {\mathcal P}^{\uparrow}_{n, k}$ are called {\it canonical}.
Similarly, permutations of the form $na_1a_2 \cdots a_{n-1} \in {\mathcal P}^{\downarrow}_{n, k}$ are also 
called {\it canonical}. Alternatively though equivalently, in view of (a)--(c), we could call those permutations of the forms
$1a_1a_2 \cdots a_{n-1}$ and $a_1a_2 \cdots a_{n-1}1$ canonical as in [4]. It is proved in [4, Theorem 6] 
that the period $\pi(\xi)$ satisfies the relation 
\begin{eqnarray}
     \pi(\xi) = (n-k)\gcd(n,\pi(\xi)),  
\end{eqnarray}
for $\xi \in {\mathcal P}^{\uparrow}_{n, k}$. \\ 
\indent
From the properties (a) and (b) of $\sigma$ we see that orbits of permutations of ${\mathcal P}^{\uparrow}_{n, k}$ 
contain canonical ones $a_1a_2 \cdots a_{n-1}n$, where $a_1a_2 \cdots a_{n-1}$ are permutations of $[n-1]$. The next 
theorem plays a fundamental role in enumerating PAPs with given ascent number. The procedure is analogous to
[5], where classical Eulerian numbers $|{\mathcal S}_{n, k}|$ are discussed.  \\
\\
{\bf Theorem 1.} {\it Let $n$ be an even integer and $k$ be such that $1 \le k \le n-1$. Then it follows that }
\begin{eqnarray*}
|{\mathcal P}^{\uparrow}_{n, k}|  =  (n-k)|{\mathcal P}_{n-1, k-1}|, ~~~|{\mathcal P}^{\downarrow}_{n,k}|
 ~= ~(k+1)|{\mathcal P}_{n-1, k}|;  \\
|{\mathcal P}_{n, k}|   =   (n-k)|{\mathcal P}_{n-1, k-1}| + (k+1)|{\mathcal P}_{n-1, k}|.~~~~~~
\end{eqnarray*}
\\
{\bf Proof.}  
Considering the orbit of $\xi \in {\mathcal P}^{\uparrow}_{n, k}$ under $\sigma$, it follows from (3) that its period is of the form
$d(n - k)$, where $d = \gcd(n,\pi(\xi))$ is a divisor of $n$. 
For a divisor $d$ of $n$, we denote by $\alpha_d^{k}$ the number of orbits with
period $d(n-k)$ in ${\mathcal P}^{\uparrow}_{n, k}$. There exist $n$ canonical permutations in
$\{ \sigma \xi, \sigma^2 \xi, \ldots, \sigma^{n(n-k)}\xi=\xi \}$ due to [4, Corollary 2], 
and hence each orbit 
\[
\{ \sigma \xi, \sigma^2 \xi, \ldots, \sigma^{d(n-k)}\xi=\xi \}
\]
of $\xi$ with period $d(n-k)$ contains exactly $d$ canonical permutations. 
This follows from the fact that the latter repeats itself $n/d$ 
times in the former. The number of all canonical permutations in ${\mathcal P}^{\uparrow}_{n,k}$
is equal to $|{\mathcal P}_{n-1, k-1}|$, since canonical permutations are represented as
$a_1a_2 \cdots a_{n-1}n$, where $a_1a_2 \cdots a_{n-1}$ are permutations of ${\mathcal P}_{n-1, k-1}$.
Since there exist $\alpha_d^{k}$ orbits with period $d(n-k)$ for each divisor $d$ of $n$,
classifying all canonical permutations of ${\mathcal P}^{\uparrow}_{n, k}$ into orbits leads us to 
\begin{eqnarray}
      |{\mathcal P}_{n-1, k-1}| =  \sum_{d|n} d \alpha_d^{k}.
\end{eqnarray}
Using the numbers of orbits and periods, we see that the cardinality of 
${\mathcal P}^{\uparrow}_{n, k}$ is given by 
\begin{eqnarray}
     |{\mathcal P}^{\uparrow}_{n, k}| =  \sum_{d|n} d(n-k) \alpha_d^{k} = (n-k)  \sum_{d|n} d\alpha_d^{k}  
     =  (n-k) |{\mathcal P}_{n-1, k-1}|.
\end{eqnarray}  
\indent
Next consider the set of all permutations $\xi = a_1a_2 \cdots a_n \in {\mathcal P}^{\downarrow}_{n, k}$, {\it i.e}., $a_1 > a_n$. 
For such $\xi = a_1a_2 \cdots a_n$ let us define the involution
$\xi^{\ast} = a_n \cdots a_2a_1$. Then the set ${\mathcal P}^{\downarrow}_{n, k}$ is converted to 
${\mathcal P}^{\uparrow}_{n, n-k-1}$ by the involution. 
Therefore, using (5), its cardinality is 
\begin{eqnarray}
    |{\mathcal P}^{\downarrow}_{n, k}| = |{\mathcal P}^{\uparrow}_{n, n-k-1}| 
    = (k+1)|{\mathcal P}_{n-1, n-k-2}| = (k+1) |{\mathcal P}_{n-1, k}|,
\end{eqnarray}  
since the involution operator $\ast$ is a bijection from ${\mathcal P}_{n-1, n-k-2}$ into ${\mathcal P}_{n-1, k}$. \\
\indent
Since all permutations $a_1a_2 \cdots a_n$ of ${\mathcal P}_{n, k}$ are divided in two classes
according to $a_1 < a_n$ or $a_1 > a_n$. Therefore, by adding (5) and (6), we obtain 
\begin{eqnarray*}
  |{\mathcal P}_{n, k}| = (n-k) |{\mathcal P}_{n-1, k-1}| + (k+1) |{\mathcal P}_{n-1, k}|.
\end{eqnarray*}
This completes the proof.   \\
\\
\indent
It is interesting to observe that the last recurrence relation holds only for even $n$, 
although such a relation is valid for classical Eulerian numbers $|{\mathcal S}_{n, k}|$ 
independently of the parity of $n$ (see [1, 3]). \\
\\
{\bf Theorem 2.} {\it Let $n$ be an even integer and $k$ be such that $1 \le k \le n-1$. Then it follows that }
\[
|{\mathcal O}^{\uparrow}_{n, k}| = |{\mathcal E}^{\uparrow}_{n, k}| + |{\mathcal P}_{n-1, k-1}|
, ~~|{\mathcal E}^{\downarrow}_{n,k}| = |{\mathcal O}^{\downarrow}_{n,k}| + |{\mathcal P}_{n-1, k}|.  
\]
\\
{\bf Proof.}  We consider the image of ${\mathcal O}^{\uparrow}_{n, k}$ under the operator $\sigma$. Permutations 
$\xi = a_1a_2 \cdots a_n$ in ${\mathcal O}^{\uparrow}_{n, k}$ begin with odd entries $a_1$ and end with
$a_n (> a_1)$.  Since $\sigma$ preserves
the ascent number, property (a) of $\sigma$ implies that the image includes all permutations of 
${\mathcal E}^{\uparrow}_{n, k}$ that begin with even entries. Moreover,
it also includes permutations of the form $1b_1b_2 \cdots b_{n-1}$, which are images of $a_1a_2 \cdots a_{n-1}n
\in {\mathcal O}^{\uparrow}_{n, k}$ by $\sigma$ (see (b) of Section 1). Then the permutation 
$(b_1-1)(b_2-1) \cdots (b_{n-1}-1)$ that is obtained by subtracting one from each entry of $b_1b_2 \cdots b_{n-1}$ 
is a PAP of $[n-1]$ and has ascents $k-1$. The number of such permutations is equal to $|{\mathcal P}_{n-1, k-1}|$. Hence we have
\[
|{\mathcal O}^{\uparrow}_{n, k}| = |{\mathcal E}^{\uparrow}_{n, k}| + |{\mathcal P}_{n-1, k-1}|.
\]
\indent
Similarly, we consider the image of ${\mathcal E}^{\downarrow}_{n, k}$ under the operator $\sigma$. Permutations 
$\xi = a_1a_2 \cdots a_n$ in ${\mathcal E}^{\downarrow}_{n, k}$ begin with even entries $a_1$ and end with
$a_n (< a_1)$.  Since $\sigma$ preserves
the ascent number, property (a) of $\sigma$ implies that the image includes all permutations of 
${\mathcal O}^{\downarrow}_{n, k}$ that begin with odd entries. Moreover,
it also includes permutations of the form $b_1b_2 \cdots b_{n-1}1$, which are images of $na_1a_2 \cdots a_{n-1}
\in {\mathcal E}^{\downarrow}_{n, k}$ by $\sigma$ (see (c) of Section 1). Then the permutation 
$(b_1-1)(b_2-1) \cdots (b_{n-1}-1)$ is a PAP of $[n-1]$ and has the same ascents $k$. 
The number of such permutations is equal to $|{\mathcal P}_{n-1, k}|$. Hence we have
\[
|{\mathcal E}^{\downarrow}_{n, k}| = |{\mathcal O}^{\downarrow}_{n, k}| + |{\mathcal P}_{n-1, k}|.
\]
This completes the proof. \\
\\
\indent
Combining Theorems 1 and 2 yields the following Corollary.\\
\\
{\bf Corollary 3.} {\it The cardinalities of the above four sets are given by}
\begin{eqnarray*}
|{\mathcal O}^{\uparrow}_{n, k}| &=& \frac{(n-k+1)|{\mathcal P}_{n-1, k-1}|}{2},  
~~ |{\mathcal E}^{\uparrow}_{n, k}| ~=~ \frac{(n-k-1)|{\mathcal P}_{n-1, k-1}|}{2},  \\
|{\mathcal E}^{\downarrow}_{n, k}| &=& \frac{(k + 2)|{\mathcal P}_{n-1, k}|}{2},  
~~ ~~~~~~|{\mathcal O}^{\downarrow}_{n, k}| ~=~ \frac{k|{\mathcal P}_{n-1, k}|}{2}.
\end{eqnarray*}
\\
{\bf Proof.} Since $|{\mathcal P}^{\uparrow}_{n, k}|  =  |{\mathcal O}^{\uparrow}_{n, k}|
 + |{\mathcal E}^{\uparrow}_{n, k}| = (n-k)|{\mathcal P}_{n-1, k-1}|$ by Theorem 1, we obtain
 \[
(n-k)|{\mathcal P}_{n-1, k-1}| = 2|{\mathcal E}^{\uparrow}_{n, k}| + |{\mathcal P}_{n-1, k-1}|
\]
by Theorem 2, which leads us to
\[
 |{\mathcal E}^{\uparrow}_{n, k}| = \frac{(n-k-1)|{\mathcal P}_{n-1, k-1}|}{2},~~
 |{\mathcal O}^{\uparrow}_{n, k}| = \frac{(n-k+1)|{\mathcal P}_{n-1, k-1}|}{2}.
\]
Similarly, Since $|{\mathcal P}^{\downarrow}_{n, k}|  =  |{\mathcal E}^{\downarrow}_{n, k}|
 + |{\mathcal O}^{\downarrow}_{n, k}| = (k+1)|{\mathcal P}_{n-1, k}|$ by Theorem 1, we obtain
\[
(k+1)|{\mathcal P}_{n-1, k}| = 2|{\mathcal O}^{\downarrow}_{n, k}| + |{\mathcal P}_{n-1, k}|
\]
by Theorem 2, which implies
\[
  |{\mathcal O}^{\downarrow}_{n, k}| = \frac{k |{\mathcal P}_{n-1, k}|}{2}, 
  ~~ |{\mathcal E}^{\downarrow}_{n, k}| = \frac{(k+2)|{\mathcal P}_{n-1, k}|}{2}.
\]
This completes the proof.
\\
\\
\indent
Just as classical Eulerian numbers $|{\mathcal S}_{n, k}|$ are divisible by some prime powers under certain conditions 
(see [5, Theorem 7]), we will show that several cardinalities related to ${\mathcal P}_{n, k}$ or
${\mathcal P}_{n-1, k}$ can be divided by some prime powers under the same conditions.\\
\\
{\bf Theorem 4.} 
{\it Suppose that $p$ is a prime and an even integer $n$ is divisible by $p^m$ 
for a positive integer $m$. If $k$ is divisible by $p$, then $|{\mathcal P}_{n-1, k-1}|$ is divisible by $p^m$ and
$|{\mathcal P}^{\uparrow}_{n, k}|$ is divisible by $p^{m+1}$. 
If $k + 1$ is divisible by $p$, then $|{\mathcal P}^{\downarrow}_{n, k}|$ is divisible by $p^{m+1}$.} \\
\\
{\bf Proof.}  We denote by $\alpha_d^{k}$ the number of orbits with period $d(n-k)$ in 
${\mathcal P}^{\uparrow}_{n,k}$, as in the proof of Theorem 1. 
Without loss of generality we can assume that $m$ is the largest integer for which $p^m$ divides $n$. 
Suppose $k$ is a multiple of $p$. From (3) the period of a permutation
$\xi \in {\mathcal P}^{\uparrow}_{n,k}$ satisfies
$\pi(\xi) = d(n-k)$, where $d=\gcd(n,\pi(\xi))$. Hence we get $d=\gcd(n,d(n-k))$ and
$\gcd(n/d, n-k) = 1$. This implies that $\gcd(k, n/d) = 1$ holds. In other words,
$\alpha_d^k = 0$ for any divisor $d$ of $n$ such that $\gcd(k, n/d) > 1$.
Moreover, a divisor $d$ with $\gcd(k, n/d) =1$ must be a multiple of $p^m$, since $k$ is a multiple of $p$. Therefore, 
equality (4) implies that $|{\mathcal P}_{n-1, k-1}|$ is divisible by $p^m$. Theorem 1 tells us that
$|{\mathcal P}^{\uparrow}_{n, k}|  =  (n-k)|{\mathcal P}_{n-1, k-1}|$. Since $n$ and $k$ are multiples of $p$ by the
assumption, $|{\mathcal P}^{\uparrow}_{n, k}|$ is divisibles by $p^{m+1}$. 
Since $|{\mathcal P}^{\downarrow}_{n, k}|  =  (k+1)|{\mathcal P}_{n-1, k}|$ holds from Theorem 1, it is also divisible by $p^{m+1}$, when
$k + 1$ is a multiple of $p$. This follows from the fact that $|{\mathcal P}_{n-1, k}|$ is divisible by $p^m$ 
when $k + 1$ is a multiple of $p$.  \\
\\
\indent
Using Corollary 3, similar arguments imply the following result. \\
\\
{\bf Corollary 5.} 
{\it Suppose that $p$ is an odd prime and an even integer $n$ is divisible by $p^m$ 
for a positive integer $m$. If $k$ is divisible by $p$, then 
$|{\mathcal O}^{\uparrow}_{n, k}|$ and $|{\mathcal E}^{\uparrow}_{n, k}|$ are divisible by $p^m$. 
If $k+1$ is divisible by $p$, then 
$|{\mathcal O}^{\downarrow}_{n, k}|$ and $|{\mathcal E}^{\downarrow}_{n, k}|$ are divisible by $p^m$.} \\
\\
\\
{\large \bf 3.  Signed Eulerian numbers} \\
\\
\indent
Let ${\mathcal P}^{\rm e}_{n,k}$ and ${\mathcal P}^{\rm o}_{n,k}$ be the sets of even permuations 
and odd ones in ${\mathcal P}_{n,k}$, respectively. The signed Eulerian number introduced in Section 1 is
defined by the difference of classical Eulerian numbers
\[  d_{n,k} = |{\mathcal S}^{\rm e}_{n,k}| - |{\mathcal S}^{\rm o}_{n,k}|.
\]
In this section we prove that they have an alternative expression as
\[  d_{n,k} =  |{\mathcal P}^{\rm e}_{n,k}| - |{\mathcal P}^{\rm o}_{n,k}|.
\]
In other words, the numbers of even permutations and odd ones are the same in the set of permutations
that are not parity alternating in ${\mathcal S}_{n,k}$.  So defining 
\[
{\mathcal N}_{n,k} = {\mathcal S}_{n,k} \setminus {\mathcal P}_{n,k}, 
\] 
the set of all permutations that are not parity alternating in ${\mathcal S}_{n,k}$, we show that  the numbers of even 
permutations and odd ones are equal in ${\mathcal N}_{n,k}$. This assertion was proved in [7] 
by using the recurrence relation (2). The present proof is self-contained and more elementary.  \\
\indent
In order to do so, we introduce another operator based on $\sigma$. Let us define an operator $\tau$ on canonical permutations of ${\mathcal S}_{n,k}$ by
\begin{eqnarray}
\tau (a_1a_2\cdots a_{n-1}n) = \sigma^{n-a_{n-1}}(a_1a_2\cdots a_{n-1}n) = b_1b_2\cdots b_{n-1}n,
\end{eqnarray}
where $b_1 =  n-a_{n-1}$ and $b_i = a_{i-1} + (n-a_{n-1})$ (mod $n$) for $i$ $(2 \le i \le n-1)$. 
Since $\sigma$ preserves the ascent number of a permutation, so does $\tau$. 
Notice that $\tau^n \xi = \xi$ holds for a canonical permutation $\xi$, since the differences between adjacent entries 
remain constant throughout, and every entry of $\xi$ returns to the original position after $n$ applications of $\tau$. \\
\indent
It is easy to see that in the case of even $n$ the operator $\tau$ maps canonical permutations of ${\mathcal P}_{n, k}$
onto themselves and similarly for ${\mathcal N}_{n, k}$.
On the other hand, when $n$ is odd, this is not true; $\xi = 1436527$, for example, is a canonical PAP, but
$\tau \xi = 5621437$ is not a PAP. So at first we assume $n$ is an even positive integer. \\
\indent
A {\it PAP sequence} $a_{i+1} \cdots a_j$ of a permutation $a_1a_2\cdots a_n$ means a consecutive subsequence 
in which even and odd integers appear alternately. A PAP sequence $a_{i+1} \cdots a_j$ is called 
{\it maximal}, if neither $a_ia_{i+1} \cdots a_j$ nor $a_{i+1} \cdots a_ja_{j+1}$ is no longer a PAP sequence. \\
\\
{\bf Lemma 6.} {\it Let $n$ be an even positive integer. Suppose that $\xi = a_1 \cdots a_ia_{i+1}  \cdots a_{n-1}n$ 
is not a PAP, while $a_{i+1}  \cdots a_{n-1}n$ is a maximal PAP sequence.
Then all of $\{\xi, \tau \xi, \ldots, \tau^{n-i-1} \xi\}$ have the same parity and 
the first entries of the permutations in $\{\tau \xi, \tau^2 \xi, \ldots, \tau^{n-i-1} \xi\}$ are all odd, but $\tau^{n-i} \xi$ has a different parity
from that of $\{\xi, \tau \xi, \ldots, \tau^{n-i-1} \xi\}$ and the first entry of $\tau^{n-i} \xi$ is even. } \\
\\
{\bf Proof.}  
First we examine the parity of $\tau \xi$. 
Remark that ${\rm inv}(\sigma \xi) = {\rm inv}(\xi)$, because $\xi$ is canonical. However,
each additional application of $\sigma$ changes the parity of permutations 
as long as $n$ lies in the interior of permutations, as shown
in Section 1. Therefore, after the entry $a_{n-1}$ of $\xi$ becomes $n$ at the
right end of  $\tau \xi$ by the application of $\sigma^{n-a_{n-1}}$, 
the parity of $\xi$ has changed $n-a_{n-1}-1$ times. Hence,
when $a_{n-1}$ is even, the parity of $\tau \xi$ is different from that of $\xi$ and 
the first entry of $\tau \xi$ is an even $n-a_{n-1}$. On the other hand,  
when $a_{n-1}$ is odd, the parity of $\tau \xi$ is the same as that of $\xi$ and 
the first entry of $\tau \xi$ is odd. \\
\indent
In the latter case, let $\tau \xi = b_1b_2 \cdots b_{n-1}n$, 
where $b_{n-1} = a_{n-2}+ (n-a_{n-1})$ (mod $n$). If $a_{n-2}$ is even, then $b_{n-1}$ is odd. So the above
argument can be applied to $\tau \xi$, implying that $\{\tau \xi, \tau^2 \xi \}$ have the same parity. 
Hence all of $\{ \xi$, $\tau \xi$, $\tau^2 \xi\}$ have the same parity. 
Moreover, the first entry of $\tau^2 \xi$ is odd. But if $a_{n-2}$ is odd, then $b_{n-1}$ is even
and $\tau^2 \xi$ has
a different parity from that of $\{\xi,  \tau\xi\}$. Moreover, 
the first entry of $\tau^2 \xi$ is even. If $a_{i+1}  \cdots a_{n-1}n$ is a maximal PAP 
sequence, we see that the lemma follows, by employing this argument $n-i$ times.  \\
\\
\indent
By Lemma 6, when $\xi$ is not a PAP, it eventually changes the parity by repeated applications of $\tau$. 
If $\xi = a_1 \cdots a_{n-1}n$ is a PAP (hence $a_1$ is automatically an odd entry), however, then
the parity of $\{\xi, \tau \xi, \ldots, \tau^{n-1} \xi\}$ remains the same and their first entries are all odd.  \\
\indent
When $\xi$ is canonical but not a PAP, it can be written as 
\begin{eqnarray}
    \xi = a_1 \cdots a_ja_{j+1} \cdots a_ia_{i+1} \cdots a_{n-1}n = A \cdot B \cdot C,
\end{eqnarray}
where $A = a_1 \cdots a_j$ and $C = a_{i+1} \cdots a_{n-1}n$ are maximal PAP sequences and
$B = a_{j+1} \cdots a_i$ is not necessarily a PAP sequence and may be empty.
$A$ (or $C$) will be called the first (or last) maximal PAP sequence of $\xi$. 
A permutation 45316278, for example, is expressed by $\xi = 45 \cdot 316 \cdot 278$. 
The length of the sequence $A$, which is equal to $j$, is denoted by
$\langle A \rangle$, similarly for the lengths of $B$ and $C$.  \\
\\
{\bf Lemma 7.} {\it In addition to the assumptions of Lemma 6, suppose that the first entry
$a_1$ of $\xi = a_1 \cdots a_ia_{i+1}  \cdots a_{n-1}n$ is even and let 
$\tau^{n-i} \xi = c_1 \cdots c_{n-i}c_{n-i+1}  \cdots c_{n-1}n$. Then $c_1 \cdots c_{n-i}$ is a maximal PAP sequence.
} \\
\\
{\bf Proof.}  
Let us put $\xi = a_1 \cdots a_ia_{i+1} \cdots a_{n-1}n=A \cdot B \cdot C$ and $\tau \xi = b_1b_2 \cdots b_{n-1}n$,
where $C = a_{i+1} \cdots a_{n-1}n$ is the last maximal PAP sequence of $\xi$.
If $a_{n-1}$ is even, or $\langle C \rangle =1$, then both $b_1 = n-a_{n-1}$ and $b_2 = a_1 + (n-a_{n-1})$ (mod $n$)
are even, because $a_1$ is even. 
If $a_{n-1}$ is odd, or $\langle C \rangle \ge 2$, then both $b_1$ and $b_2$ are odd. This implies that 
the first PAP sequence of $\tau \xi$ necessarily has a length of one. \\
\indent
By each application of $\tau$ to $\xi$, the entries 
$\{n, a_{n-1}, \ldots, a_{i+1}\}$ of $C$ move to the left end of permutations one by one
in this order.
If we put $\tau^{n-i} \xi = c_1 \cdots c_{n-i}c_{n-i+1} \cdots c_{n-1}n$, then 
$c_1 \cdots c_{n-i}$ is a PAP sequence and it is maximal. 
This follows from the facts that $a_{i+1}  \cdots a_{n-1}n$ is a PAP sequence of $\xi$ and that
the parity of $c_{n-i}$ and $c_{n-i+1}$ is the same, because $\tau^{n-i} \xi = \tau^{n-i-1}(\tau \xi)$ holds and both
$b_1$ and $b_2$ of $\tau \xi$ have the same parity. \\
\\
\indent
Lemma 7, together with Lemma 6, states that applying $\tau^{n-i}$ to $\xi$ of (8) moves
the last maximal PAP sequence, $C$, of $\xi$ to the first maximal PAP sequence
of a permutation $\tau^{n-i} \xi$ with opposite parity and both sequences have the same length. 
Furthermore, all entries of $A$ and $B$ move to the right by
$n-i$ positions by means of an application of $\tau^{n-i}$, although their values have changed. 
In particular, the last entry of $B$, $a_i$, turns to $n$ at the right end of a permutation,
when $B$ is not empty.  \\
\indent
As an illustrative example, applying $\tau^3$ to an odd permutation $\xi = 45 \cdot 316 \cdot 278$ yields an even permutation
$412 \cdot 675 \cdot 38$. Thus the last maximal PAP sequence $278$ of $\xi$ moves to the first maximal one $412$ of 
$\tau^3 \xi$, and the remaining sequence $45316$ of $\xi$ moves to the right as $67538$ in $\tau^3 \xi$. 
\\
\\
{\bf Theorem 8.} {\it In ${\mathcal N}_{n,k}$ the total number of even permutations is equal to that of odd ones.} \\
\\
{\bf Proof.}  
Part I (the case of odd $n$). \\
\indent
For odd $n$ let us introduce the set of all canonical permutations in ${\mathcal N}_{n+1, k+1}$,
which we denote by ${\mathcal N}^{\rm c}_{n+1, k+1}$ and on which we will apply the operator $\tau$ of (7). 
Then it is possible to correspond  
$\bar{\xi} =a_1a_2\cdots a_n(n+1)$ in ${\mathcal N}^{\rm c}_{n+1, k+1}$ 
to each $\xi = a_1a_2\cdots a_n$ in ${\mathcal N}_{n, k}$.
Note that ${\rm inv}(\xi) = {\rm inv}(\bar{\xi})$ and the correspondence $\xi \rightarrow \bar{\xi}$
is a bijection from ${\mathcal N}_{n, k}$ to ${\mathcal N}^{\rm c}_{n+1, k+1}$. Therefore, the number of even 
permutations in the latter is equal to that in the former and similarly for odd permutations.
Since $n+1$ is even, we can apply Lemmas 6 and 7 to canonical permutations of ${\mathcal N}^{\rm c}_{n+1, k+1}$, 
in which it is possible to find those beginning with even $a_1$. 
We denote by ${\mathcal N}^{\rm c \dagger}_{n+1, k+1}$ the set of such canonical ones having the
even first entry. \\
\indent
A permutation $\bar{\xi} = a_1a_2\cdots a_n(n+1)$ in ${\mathcal N}^{\rm c}_{n+1, k+1}$ can be written as
$\bar{\xi} = A \cdot B\cdot C$ by using maximal PAP sequences $A = a_1\cdots a_j$ and 
$C = a_{i+1} \cdots (n+1)$, as in (8). 
So we have $n+1 = \langle A \rangle+\langle B \rangle+\langle C \rangle$. 
The sequence $B$ needs not be a PAP sequence and may be empty. As was proved in Lemma 6, 
until the operator $\tau$ is applied $\langle C \rangle$ times to $\bar{\xi}$, the parity does not change. \\
\indent
Suppose that $\bar{\xi}_1, \bar{\xi}_2, \ldots, \bar{\xi}_m$ are all of even canonical permutations in 
${\mathcal N}^{\rm c \dagger}_{n+1, k+1}$, {\it i.e.}, even permutations beginning with even integers and ending with $n+1$.
For each $i~(1 \le i \le m)$ let us define $\ell_i$ as the smallest positive integer such that 
$\bar{\eta}_i = \tau^{\ell_i}\bar{\xi}_i$ becomes an odd permutation in ${\mathcal N}^{\rm c}_{n+1, k+1}$, 
which indeed belongs to ${\mathcal N}^{\rm c \dagger}_{n+1, k+1}$ by Lemma 6.
The correspondence between $\bar{\xi}_i$ and $\bar{\eta}_i$ is a bijection and thus we get all of odd ones 
$\bar{\eta}_1, \bar{\eta}_2, \ldots, \bar{\eta}_m$ in ${\mathcal N}^{\rm c \dagger}_{n+1, k+1}$ by this procedure.  
For $i~(1 \le i \le m)$ let $\bar{\xi}_i = A_i\cdot B_i\cdot C_i$  and 
$\bar{\eta}_i = A_i^{\natural}\cdot B_i^{\natural}\cdot C_i^{\natural}$ be the
expressions using maximal first and last PAP sequences. Then it follows that the PAP sequence $C_i$ of 
$\bar{\xi}_i$ moves into the PAP sequence $A_i^{\natural}$ of $\bar{\eta}_i$ by $\tau^{\ell_i}$, and 
$\ell_i = \langle C_i \rangle = \langle A_i^{\natural} \rangle$ from Lemma 7. 
So we have
\begin{eqnarray}
        \sum_{i=1}^m \langle A_i^{\natural} \rangle = \sum_{i=1}^m \langle C_i \rangle.
\end{eqnarray}
\indent
Since each odd $\bar{\eta}_i=A_i^{\natural}\cdot B_i^{\natural}\cdot C_i^{\natural}$ is, in turn,  
changed into a certain even permutation 
$\bar{\xi}_j = A_j \cdot B_j \cdot C_j$ in ${\mathcal N}^{\rm c \dagger}_{n+1, k+1}$ in a similar way, we obtain 
an analogous equality
\begin{eqnarray}
        \sum_{i=1}^m \langle A_i \rangle = \sum_{i=1}^m \langle C_i^{\natural} \rangle.
\end{eqnarray}
Moreover, the PAP sequence $C_i^{\natural}$ of $\bar{\eta}_i$ moves into the PAP sequence
$A_j$ of $\bar{\xi}_j$. From the above equalities we also have
\begin{eqnarray*}
        \sum_{i=1}^m {\langle B_i \rangle} = \sum_{i=1}^m {\langle B_i^{\natural} \rangle}.
\end{eqnarray*}	
For each $i~(1 \le i \le m)$ the permutation $\tau^{\ell}\bar{\xi}_i$ remains even, whenever 
$\ell < \langle C_i \rangle$.
Noting $\tau^0\bar{\xi}_i = \bar{\xi}_i$, the number of all even permutations in ${\mathcal N}^{\rm c}_{n+1, k+1}$ is 
given by the sum (9). 
Similarly, the number of all odd permutations in ${\mathcal N}^{\rm c}_{n+1, k+1}$ is given by the sum (10). \\
\indent
Now we know that the last maximal PAP sequence $C_i$ of $\bar{\xi}_i$ moves into the first maximal 
PAP sequence $A_i^{\natural}$ of $\bar{\eta}_i$ and 
that the last maximal PAP sequence $C_i^{\natural}$ of $\bar{\eta}_i$ moves into the first maximal PAP
sequence $A_j$ of some even $\bar{\xi}_j$.
Thus we obtain
\[
\sum_{i=1}^m \langle C_i \rangle=\sum_{i=1}^m \langle C_i^{\natural} \rangle,
\]
and hence the values of (9) and (10) are equal; $\sum_{i=1}^m \langle A_i \rangle=\sum_{i=1}^m \langle C_i \rangle$.  \\
\indent
This implies that, among ${\mathcal N}^{\rm c}_{n+1, k+1}$, the number of even permutations 
is equal to that of odd ones.
In other words, the number of even permutations in ${\mathcal N}_{n, k}$ is equal to that of odd ones in it, which
completes the proof of Part I. \\
\\
Part II (the case of even $n$). \\
\indent
The case of even $n$ can easily be proved from the former part. Notice that
in this case $\xi \in {\mathcal N}_{n, k}$ if and only if $\sigma\xi \in {\mathcal N}_{n, k}$.
We divide all even permutations $a_1a_2 \cdots a_n$ of ${\mathcal N}_{n, k}$ 
into the three types by the position of $n$: 
\begin{itemize}
\item[(i)]  $a_i = n$ for some $i$ $(2 \le i \le n-1)$;
\item[(ii)]  $a_n = n$;
\item[(iii)]  $a_1 = n$.
\end{itemize}
\indent
\indent
On the other hand, we divide all odd permutations $a_1a_2 \cdots a_n$ of 
${\mathcal N}_{n, k}$ into the following three types by the position of one: 
{\begin{itemize}
\item[(i')]  $a_i = 1$ for some $i$ $(2 \le i \le n-1)$;
\item[(ii')]  $a_1 = 1$;
\item[(iii')]  $a_n = 1$. 
\end{itemize}}
\indent
Let $\xi = a_1a_2 \cdots a_n$ be an even permutation of type (i). 
Then, using property (a) of $\sigma$, we see that $\sigma \xi$ is an odd one of type (i'), because 
the difference of the numbers of inversions between $\xi$ and $\sigma \xi$ is
$|2i- (n+1)|$ by (1) and it is odd by assumption. Since $\sigma$ is a bijection, we see that to each even
permutation $\xi$ of type (i) in ${\mathcal N}_{n, k}$ corresponds an odd one 
of type (i') in ${\mathcal N}_{n, k}$. Therefore, both types have the same cardinality.  \\
\indent
Let $\xi = a_1a_2 \cdots a_{n-1}n$ be of type (ii), where $a_1a_2 \cdots a_{n-1}$ is a
permutation of  $[n-1]$. Since $\xi$ is an even permutation, 
the permutation $a_1a_2 \cdots a_{n-1} \in {\mathcal N}_{n-1, k-1}$ is also even. 
Hence the number of elements in ${\mathcal N}(n,k)$ of type (ii) is the cardinality
of all even permutations of ${\mathcal N}(n-1,k-1)$. On the other hand, the set ${\mathcal N}_{n, k}$ of type (ii')
consists of all odd $\xi = 1a_2a_3 \cdots a_n$ and hence the permutations $(a_2-1)(a_3-1) \cdots (a_n-1)$, 
that are obtained by subtracting one from the entries $a_2a_3 \cdots a_n$, are all odd in 
${\mathcal N}_{n-1, k-1}$. Using Part I, we see that both cardinalities are the same. \\ 
\indent
Similar arguments can be applied to types (iii) and (iii'). 
Let $\xi = na_1a_2 \cdots a_{n-1}$ be of type (iii). 
Since $n-1$ is odd and $\xi$ is an even permutation, we see that the permutation $a_1a_2 \cdots a_{n-1}$ is odd in 
${\mathcal N}_{n-1, k}$. 
Hence the number of even permutations in ${\mathcal N}_{n, k}$ of type (iii) is equal to that 
of all odd permutations of ${\mathcal N}_{n-1, k}$.
On the other hand, the set ${\mathcal N}_{n, k}$ of type (iii')
consists of all odd $\xi = a_2a_3 \cdots a_n1$, where $(a_2-1)(a_3-1) \cdots (a_n-1)$ are all even in 
${\mathcal N}_{n-1, k}$. Hence it follows from Part I that both cardinalities are the same. 
Thus we see that the cardinality of even permutations in ${\mathcal N}_{n, k}$ is equal to that of odd ones, which  
completes the proof. \\
\\
\indent
Finally we leave an unsolved problem, which is a counterpart of Theorem 1; is it possible to express $|{\mathcal P}_{n, k}|$ for odd $n$ by means
of $|{\mathcal P}_{n-1, k}|$ {\it etc.} in a relatively simple manner? \\
\\
\\
{\large \bf References} 
\begin{itemize}
\item[{[1]}]  M. Aigner, {\it Combinatorial Theory}, Springer-Verlag, 1979.
\item[{[2]}]  J. D\'esarm\'enien and D. Foata, The signed Eulerian numbers, {\em Discrete Mathematics}
      {\bf 99}, 49--58, 1992.      
\item[{[3]}]  J. Riordan, {\it An Introduction to Combinatiorial Analysis}, Princeton Univ. Press, 1980.
\item[{[4]}]  S. Tanimoto,  An operator on permutations and its 
  application to Eulerian numbers, {\em European Journal of Combinatorics} {\bf 22}, 569--576, 2001. 
\item[{[5]}]  S. Tanimoto,  A study of Eulerian numbers by means of an operator on permutations,
   {\em European Journal of Combinatorics} {\bf 24}, 33--43, 2003.
\item[{[6]}] S. Tanimoto, A study of Eulerian numbers for permutations in the alternating group, 
  {\em Integers} {\bf 6}, A31, 2006.
\item[{[7]}] S. Tanimoto, Parity alternating permutations and signed Eulerian numbers, 
  {\em Annals of Combinatorics} {\bf 14}, 355--366, 2010.
\end{itemize}
\end{document}